\documentclass[11pt]{article}
\usepackage{amsmath,amstext,amssymb,amsthm, amsfonts}
\usepackage[title]{appendix}
\usepackage{color}
\usepackage{mathrsfs}
\newtheorem{theorem}{Theorem}[section]
\newtheorem{lemma}[theorem]{Lemma}

\newtheorem{definition}[theorem]{Definition}
\newtheorem{assumption}[theorem]{Assumption}
\newtheorem{example}[theorem]{Example}

\newtheorem{corollary}[theorem]{Corollary}

\theoremstyle{remark}

\numberwithin{equation}{section}

\def\P{\mathscr{P}} 		
\def\B{\mathcal{B}}

\def\E{\mathbf{E}}
\def\X{\mathbb X}

\def\C{\mathcal}
\def\BB{\mathbf}
\def\F{\mathfrak}

\def\F{\mathcal{F}}

\begin{document}

\title{Kolmogorov's Equations for Jump Markov Processes
and their Applications to Control Problems}

\date{}

\maketitle

\begin{center}
Eugene~A.~Feinberg\footnote{Department of Applied Mathematics and
Statistics,
 Stony Brook University,
Stony Brook, NY 11794-3600, USA, eugene.feinberg@sunysb.edu},
Albert~N.~Shiryaev\footnote{Steklov Mathematical Institute and Moscow State University,
Moscow, Russia, albertsh@mi-ras.ru}
\end{center}

\begin{abstract}
This paper describes the structure of solutions to Kolmogorov's equations for nonhomogeneous jump Markov processes and applications of these results to control of jump stochastic systems.  These equations were studied by Feller (1940), who clarified in 1945 in the errata to that paper that some of its results covered only nonexplosive Markov processes. In this work, which is largely of a survey nature,
the case of explosive processes is also considered.  This paper is based on the invited talk presented by the authors at the conference ``Chebyshev-200", and it describes the results of their joined studies with Manasa Mandava (1984-2019).
\end{abstract}

Key words: Kolmogorov's equations, jump Markov processes, optimal control.

\maketitle


\section{Introduction}\label{s1}

\textbf{1.1.} A. N. Kolmogorov introduced backward and forward equation in his seminal paper  \cite{Kol} ``On analytic methods in probability theory"  published in 1931.  Among many deep ideas and results introduced there, he explicitly wrote backward and forward equations for jump Markov processes with finite and countable state spaces and for diffusion processes.

W. Feller studied Kolmogorov's equations in several publications. In  partucular, in \cite{Fel} he studied Kolmogprov's equations for nonhomogeneous jump Markov processes with Polish state spaces. He clarified later in the errata to \cite{Fel} that the uniqueness results in \cite{Fel} were correct only for nonexplosive processes. Examples of multiple solutions to Kolmogorov's equations for homogeneous jump Markov processes with countable state spaces were provided by Doob~\cite{Do}, Kendall~\cite{Ken}, and Reuter~\cite{Reu}; see also the book by Anderson~\cite{Anderson}.

\textbf{1.2.}  For possibly explosive processes
  Feller~\cite{Fel} wrote explicit formulae for transition probabilities.  In this paper we use the term ``transition function" instead of a transition probability because the full measure may not be equal to 1.

  It is easy to conclude from Feller's~\cite{Fel} arguments that the transition function of the jump Markov process is the \emph{minimal} solution of the corresponding Kolmogorov \emph{backward} equation.  Feller~\cite{Fel} mostly studied backward equations trying to show that a solution of a backward equation also solves a forward equation for sets of states with bounded jump intensities.  The possible reason for this approach, as this reason was mentioned by Feller~\cite{Fel}, was that, while Kolmogorov's backward equations can be written under very general conditions, the right-hand side of Kolmogorov's forward equations can be an uncertainty of the form of $\infty-\infty$ for some sets of states even if jump rates are bounded at each state.

\textbf{1.3.} There are many applications of jump Markov processes, and one of them is the area of continuous-time jump Markov decision processes (CTJMDPs); see monographs \cite{GH, KR, PZ} dealing with this topic. One of the basic questions for CTJMDPs is how to define a policy and a stochastic process defined by a policy.

The first publications on CTJMDPs \cite{Ho,Ka,Mi,Mi1,Ry} dealt either with stationary or Markov policies, that is, decisions depend either only on the current state or on the current state and time, and stochastic processes, which in these cases are jump Markov, were defined in these papers via solutions of Kolmogorov's forward equations.

Feller's~\cite{Fel} results were used in some studies including \cite{GHL,GR} to define jump Markov processes for CTJMDPs with Polish state spaces.  Since Feller~\cite{Fel} used the assumption that jump intensities are continuous in the time parameter, \cite{GHL,GR} and several other publications used the assumption that Markov policies could use only decisions \emph{continuous} in the time parameter.

\textbf{1.4.}  Yushkevich~\cite{Yus, Yus1} introduced general non-Markovian policies and constructed the corresponding stochastic processes by using the Ionescu Tulcea theorem.  Kitaev~\cite{Ki} described an equivalent construction of general policies by using Jacod's~\cite{Jac} results on dual predictable projections (also called compensators) of random measureas for multivariate point processes. Kitaev~\cite{Ki} also introduced an analog of Kolmogorov's forward equations for CTJMDPs controlled by general policies; see Lemma~\ref{l:GKE} below.

In many cases, CTJMDPs can be reduced to discrete-time Markov decision processes (MDPs).  One of such reduction schemes is based on the uniformization \cite{Se}.  Another reduction scheme \cite{Fe1, Fe2} is based on the property of nonstationary exponential distributions \cite{Fe}, which implies that for infinite-horizon problems it is possible to change decisions only at jump epochs. The second reduction scheme works only when jump rates are greater than a positive number, as this always holds for problems with discounting because the discount rate can be interpreted as a jump rate.  To deal with the situations, when jump rates can be close to 0 or equal to 0, Piunovskiy~\cite{Pi} introduced an additional artificial Poisson process, and Xin Guo and Zhang~\cite{GZ} recently addressed this issue in an elegant way by doubling the state space.

\textbf{1.4.} One of the basic facts for discrete-time MDPs is that for a given initial state distribution for every policy it is possible to construct a Markov policy such that the corresponding stochastic sequences of state-action pairs for the original and new policies have the same marginal distributions \cite{DS, St}. The new Markov policy chooses the same distribution of actions as the original policy would select under the condition that only the current time, state, and the initial state distribution are known.  This fact shows that given an initial stated distribution, Markov policies are as good as general ones for basic performance criteria, and this significantly simplifies the theory of MDPs. A similar Markov policy can be also constructed for CTJMDPs, and the natural question is whether the new policy defines a stochastic process with the same marginal distributions of state-action pairs as for the original policy.  For discrete time the proofs are based on induction, and this method is not applicable for continuous time.  However, for continuous time marginal distributions of states for the original and new policies satisfy the same Kolmogorov's equations (see Lemma~\ref{solutions}(i) below), and this is the reason we started to study Kolmogorov's equations for jump Markov processes.

 In \cite{FMS} we described the structure of solutions of Kolmogorov's equations for problems with Polish state spaces approached by Feller~\cite{Fel}.  We showed in \cite{FMS} that transition functions  of jump Markov processes are the minimal solutions of Kolmogorov's backward and forward equations.  For problems with countable state spaces these results were established in Ye et al.~\cite{YGHL}.  As we mentioned above, the additional complication in the case of uncountable state spaces is that right-hand sides of Kolmogorov's forward equations can be undefined for some measurable sets of states.  We showed in \cite{FMS} that the values of the right-hand side for certain sets of states, for which the right-hand side of Kolmogorov's forward equation is well-defined, completely determine the solutions of Kolmogorov's forward equations.  In \cite{FMS2} we extended the results from \cite{FMS} to more general transition intensities; see Assumptions~\ref{Feller}-\ref{L1} below.  In \cite{FMS3} we studied applications to CTJMDPs.

Section~\ref{s2} of this paper provides main definitions and assumptions for jump Markov processes and describes the construction of transition functions.  Section~\ref{s3} describes the structures of solutions to Kolmogorov's backwards and forward equations.  The results for Kolmogorov's backward equations are currently known under more general assumptions than the similar results for Kolmogorov's forward equations.  Section~\ref{s4} deals with CTJMDPs, and it shows that the equality for marginal distributions analogous to the one that holds in discrete time, also holds for CTJMDPs, if the jump Markov process, defined by the corresponding Markov policy, is non-explosive.  If this process is explosive, then  the values of marginal probabilities of states are not greater for the corresponding Markov policies than for the original ones.  Section~\ref{s5} deals with applications of these results to particular performance criteria.

\section{ $Q$-Functions and Jump Markov Processes}\label{s2} 

\textbf{2.1.} We consider stochastic processes with values in a standard Borel space $(\X, \B(\X))$ (called the state space)  defined on the time interval $[T_0, T_1),$ where $T_0$ is a real number and $T_0<T_1\le +\infty.$ In particular, it is possible that $[T_0, T_1)=\BB{R}_+:= [0, \infty)$.
A function $P(u, x; t, B)$ with values in $[0,1],$ where $u \in [T_0, T_1)$, $t \in (u, T_1)$, $x \in \X,$ and $B \in \B(\X),$ is called  a
transition function, if it satisfies the following properties:
\begin{itemize}
\item[(i)] for all $u,x,t$ the function $P(u,x;t,\cdot)$ is a
measure on $(\X, \mathcal{B}(\X))$;
\item[(ii)]for all $B$ the
function $P(u,x; t, B)$ is Borel measurable in $(u,x,t);$
\item[(iii)]
$P(u,x;t,B)$ satisfies the Chapman-Kolmogorov equation
\begin{equation}
\label{CKE} P(u, x; t, B) = \int_{\X} P(s,y; t, B)P(u, x; s, dy),
\qquad u < s < t.
\end{equation}
\end{itemize}
A transition function $P$ is called  {\it
regular} if $P(u,x;t,\X) = 1$ for all $u,x,t$ in the domain of $P$.

A stochastic process $\xi_t: t \in [T_0, T_1)\}$ with values in $\X$,
defined on a filtered probability space $(\Omega, \C{F}, (\C{F})_{t\in[T_0,T_1)}, {\mathbf P}),$  where $(\C{F})_{t\in[T_0,T_1)}$ is a nondecreasing right-continuous family of $\sigma$-subalgebras in $\C{F}.$ is called a \emph{Markov process}  if
\begin{equation}\label{eq2.2} {\mathbf P}(\xi_t \in B \mid \C{F}_u) = {\mathbf P}(\xi_t \in B \mid
\xi_u)\quad \text{${\mathbf P}$-a.s.}\end{equation} for all $u \in [T_0, T_1),$ $t \in (u,T_1)$, and  $B
\in \B(\X)$.  As shown by Kuznetsov~\cite{Kuz}, each Markov process has a transition function $P$ such
that \[{\mathbf P}(\xi_t \in B \mid \xi_u) = P(u,\xi_u; t, B)\quad
\text{${\mathbf P}$-a.s.}\] 
This fact establishes the equivalence of two
definitions of a Markov process --- the standard definition via the formula \eqref{eq2.2} and as a process, whose transition probability satisfies the Chapman-Kolmogorov equation \eqref{CKE}. We notice that Kolmogorov~\cite{Kol} used the term \emph{stochastically defined process} instead of the term \emph{Markov process.}

\textbf{2.2.} We recall  that a \emph{multivariate point process} on a measurable filtered space $(\Omega, \mathcal{F},\mathcal{F}_{t\in(T_0,T_1)}),$ where   $(\C{F}_t)_{t \in (T_0,T_1)}$ is a nondecreasing right-continuous family of $\sigma$-subalgebras in $\mathcal{F},$ is a \emph{stochastic sequence} $(t_n,x_n)_{n=1,2,\ldots},$ where $t_n\in (T_0,T_1]$ are  stopping times
$(\{t_n \le t\} \in\mathcal{F}_t,$ $t \in (T_0, T_1]),$  $x_n\in \X\cup\{x_\infty\}$ are $\C{F}_{t_n}$-measurable  with $x_\infty\notin\X$  being an isolated point, and the following two properties hold: (a) if $t_n<T_1,$ then $t_n<t_{n+1} $ and $x_n\in \X,$ and (b) if $t_n=T_1,$ then $x_{n}=x_\infty$ and $t_{n+1}=T_1,$ where $n=1,2,\ldots;$ see \cite{JS}.

A stochastic process $\{\xi_t,\ t \in [T_0, T_1)\}$ is called a \emph{jump process} if $\xi_t= x_{n-1}$ for $t\in [t_{n-1},t_n),$ and $\xi_t= x_\infty$ for $t\ge t_\infty,$ where $t_0=T_0,$ $x_0\in\X$ is an $\C{F}_{T_0}$-measurable random variable, and the sequence $(t_n,x_n)_{n=1,2,\ldots}$ is a multivariate point process, and $t_\infty:=\lim_{n\to\infty}t_n.$ 
By defdinition,
a \emph{jump Markov process} is a Markov process which is also a jump process.

\textbf{2.3.} A function $q(x,t,B)$, where $x \in \X$, $t \in [T_0, T_1)$, and $B \in \B(\X)$,
is called a {\it Q-function} if it satisfies the following
properties:
\begin{itemize}
\item[(a)]for all $x,t$ the function $q(x,t,\cdot)$ is a signed
measure on $(\X, \mathcal{B}(\X))$ such that $q(x,t,\X)$ $\leq$ $0$ and
$0 \leq q(x, t, B \setminus \{x\}) < \infty$ for all $B \in
\B(\X)$;
\item[(b)] for all $B$ the
function $q(x,t,B)$ is measurable in $(x,t).$
\end{itemize}
If $q(x,t,\X) =0$ for all
$x\in \X$ and $t\in [T_0, T_1)$ in addition to properties (a) and (b),  then the $Q$-function $q$ is  called {\it conservative}.
Note that any $Q$-function can be transformed into a conservative
$Q$-function  by adding an additional isolated absorbing state $x'$ to $\X,$ considering the new state space $(\X',\mathcal{B}(\X'))$ with $\mathcal{B}(\X')=\sigma(\BB(\X),\{x'\}), $ and by setting
$q(x,t,E):= q(x,t,E\setminus\{x'\})-q(x,t,\X)I\{x'\in E\}$ and $q(x', t, E) := 0$ for all $E\in \mathcal{B}(\X'),$ where $x \in \X$ and $t \in
[T_0, T_1)$. Additional arguments are provided in \cite[Remark 4.1]{FMS}.

 To simplify the presentation, in this paper we
always assume that $q$ is conservative. A $Q$-function $q$ is called {\it continuous}
if it is continuous in $t \in [T_0, T_1)$.  We remark that different authors use the term $Q$-function for different objects and properties; e.g.,   this term has a different meaning in   \cite{Anderson}.

\textbf{2.3.} Let $q(x,t):=q(x,t,\X\setminus\{x\})$ for $x\in\X$ and $t\in [T_0,T_1),$ and $\bar{q}(x):\sup_{t\in [T_0,T_1)} q(x,t).$ Let us consider the following assumptions on a $Q$-function $q.$
\begin{assumption}[\cite{Fel}]
\label{Feller} For all $n = 1,2,\ldots$
there exist Borel subsets $B_n\in\mathcal{B}(\X)$ such that: (i) $B_n\subset B_{n+1}$, (ii) $B_n \uparrow \X$ as $n \to \infty,$ and  (iii) $\sup_{x \in B_n} {\bar q}(x) < n.$
\end{assumption}
\begin{assumption}[boundedness of $q(x,\cdot)$]
\label{LB}
${\bar q}(x) < \infty $ for each $x \in \X$.
\end{assumption}
\begin{assumption}[local boundedness of $q(x,\cdot)$]
\label{ALB}
$\sup_{t\in [T_0,s)} q(x,t)<\infty$ for each $s\in (T_0,T_1)$ and  $x\in \X.$
\end{assumption}

\begin{assumption}[local $\mathcal{L}^1$ boundedness of $q(x,\cdot)$]
\label{L1}
$\int_{T_0}^s q(x,t) dt < \infty$  for each $s \in (T_0, T_1)$ and  $x \in \X.$
\end{assumption}

The following lemma compares Assumptions~\ref{Feller}--\ref{L1}.
\begin{lemma}[\cite{FMS2}]
\label{l:A-eq} The following statements hold:

(i) Assumptions~\ref{Feller} and \ref{LB} are equivalent;

(ii) Assumption~\ref{LB} implies Assumption~\ref{ALB};

(iii) Assumption~\ref{ALB} implies Assumption~\ref{L1}.
\end{lemma}
\begin{proof}
 The only nontrivial statement in the lemma is that Assumption~\ref{LB} implies Assumption~\ref{Feller}, and we provide its proof. In view of \cite[Proposition 7.47]{BS}, the function $\bar{q}(x)$  is upper semianalytic.  This means that the sets $A_n:=\{x\in\X:\bar{q}(x)\ge n\}$ are analytic.   Assumption~\ref{LB} implies that
$\cap_{n=1}^\infty A_n=\emptyset.$ Novikov's separation theorem \cite[p. 219]{Kechris} implies the existence of Borel sets $C_n\supset A_n$ such that $\cap_{n=1}^\infty C_n=\emptyset.$  Let $B_n:=\cup_{m=1}^n C_m^c,$ where $D^c$ is a complement of the set $D.$ Then $B_n\uparrow\X$ and $\sup_{x\in B_n} \bar{q}(x)<n$ for all $n.$ \end{proof}

\textbf{2.5.} Let $x_\infty\notin \X$ be the point described above, and let ${\bar \X}=\X\cup\{x_\infty\}$. 
Consider the Borel $\sigma$-field $\B({\bar \X})=\sigma(\mathcal{B}(\X),\{x_\infty\})$ on $\bar \X$, which is
the minimal $\sigma$-field containing $ \mathcal{B}(\X)$ and
$\{x_\infty\}.$ Let $({\bar \X} \times (T_0, T_1])^\infty$ be the set of all sequences $(x_0, t_1, x_1,
t_2, x_2, \ldots)$ with $x_n\in \bar{\X}$ and $t_{n+1}\in
(T_0, T_1]$ for all $n =0,1,\ldots\ .$ This set is endowed
with the $\sigma$-field generated by the products of the Borel
	$\sigma$-fields $\B(\bar{\X})$ and $\B((T_0, T_1])$.

We denote by $\Omega$ the subset of all sequences $\omega= (x_0, t_1,
x_1, t_2, x_2, \ldots)$ from $({\bar \X} \times (T_0, T_1])^\infty$ such that: (i) $x_0 \in \X$; (ii) for all $n =1,2,\ldots\,,$ if $t_n < T_1$, then $t_n < t_{n+1}$ and
$x_n \in \X$, and if $t_n = T_1$, then $t_{n+1} = t_n$ and
$x_n = x_\infty$. Observe that $\Omega$ is a measurable subset of
$({\bar \X} \times (T_0, T_1])^\infty$. Consider the measurable space $(\Omega,
\C{F})$, where $\C{F}$ is the $\sigma$-field of the measurable subsets
of $\Omega$. For all $n = 0,1,\ldots$, let $x_n(\omega)=x_n$ and
$t_{n+1}(\omega)=t_{n+1},$ where $\omega \in \Omega,$ be the random variables defined on the
measurable space $(\Omega, \C{F})$. Let again $t_0 := T_0$, and  $t_\infty
(\omega) := \lim\limits_{n \to \infty} t_n (\omega)$, $\omega \in \Omega.$ For all $t \in [T_0, T_1],$ let $\C{F}_t := \sigma(\B(\X), \C{G}_t)$, where $\C{G}_t := \sigma (I\{x_n \in B\}I\{t_n \le s\}: n \ge 1, T_0 \le s \le t, B \in \B(\X)).$ Throughout this paper, we omit $\omega$ whenever possible.

Consider the multivariate point process $(t_n, x_n)_{n =1,2,\ldots}$ on $(\Omega, \C{F})$. For a given iven a $Q$-function $q$ satisfying Assumption~\ref{L1}, define a random measure $\nu$ on $([T_0, T_1) \times \X)$ by
\begin{equation}
\label{compensator}
\nu(\omega; [T_0,t], B): = \int_{T_0}^{t}\sum_{n \ge 0}I\{t_n < s \leq t_{n+1}\} q(x_n,s, B \setminus \{x_n\})ds,  \quad  t \in [T_0, T_1),\ B \in \B(\X),
\end{equation}
where $t_0=T_0$ and $\omega=x_0,t_1,x_1,t_2,\ldots\in\Omega.$
As explained in detail in \cite{FMS}, the random measure $\nu$ is predictable.
 Furthermore, $\nu(\{t\}\times \X) \le 1$ for all $t\in (T_0,T_1)$ and $\nu([t_\infty,\infty)\times \X)=0. $  According to Jacod~\cite[Theorem~3.6]{Jac}, the predictable random measure $\nu$ defined in \eqref{compensator} and a  probability measure $\gamma$ on $\X$ define a unique probability measure ${\mathbf P}$ on $(\Omega, \C{F})$ such that ${\mathbf P}(x_0 \in B) = \gamma(B), B \in \B(\X),$ and $\nu$ is the compensator of the random measure of the multivariate point process $(t_n,x_n)_{n=1,2\ldots}$ defined by the triplet $(\Omega, \C{F}, {\mathbf P})$.

Consider the process $\{\xi_t: t\in [T_0, T_1)\}$,
\begin{equation}
\label{jump}
\xi_t(\omega) := \begin{cases} x_n, &{\rm  if}\ t\in [t_n,t_{n+1})\  {\rm for}\ n=0,1,\ldots,\\ x_\infty,  & {\rm if}\  t\ge t_\infty,\end{cases}
\end{equation}
defined on $(\Omega, \C{F}, {\mathbf P})$ and adapted to the filtration $(\C{F}_{t})_{t \in [T_0, T_1)}$. By definition, the process $\{\xi_t:	 t \in [T_0, T_1)\}$ is a \emph{jump} process.

Following Feller \cite[Theorem 2]{Fel}, for $x \in \X$, $u \in [T_0, T_1)$, $t \in (u,T_1)$, and $B \in \B(\X)$, we define
\begin{equation}
\label{b0} \overline{P}^{(0)} (u,x;t,B): = I\{x \in B\} e^{-\int_u^t
q(x, s) ds}
\end{equation}
and 
\begin{equation}
\label{bn}
\overline{P}^{(n)}(u, x; t, B): = \int_{u}^{t} \int_{\X  }
 e^{ -\int_{u}^{w} q(x,\theta) d\theta }  q(x,w,dy \setminus \{x\}) \overline{P}^{(n-1)}(w, y; t, B)
 dw, \  n=1,2,\ldots,
\end{equation}
where $w \in [T_0, T_1)$ we use the notation
 $q(x,w,dz\setminus\{x\}):=q^+(x,w,dz)$ for the measure $q^+(x,w,B):= q(x,w,B\setminus \{x\})$   on the space $(\X, \B(\X))$ where  $B \in \B(\X).$
Then
\begin{equation}
\label{def}
\overline{P}(u, x; t, B) := \sum\limits_{n=0}^{\infty} \overline{P}^{(n)}(u, x; t, B)
\end{equation}
is a transition function if the $Q$-function $q$ satisfies Assumption~\ref{L1}; see \cite[(5)-(8)]{FMS2} for details, where most of the arguments are taken from \cite{Fel}. In particular, equation \eqref{bn} can be rewritten as  
\begin{equation}
\label{bn-alt}
\overline{P}^{(n)}(u,x;t,B) = \int\limits_u^t \int\limits_\X \int\limits_{B } e^{-\int_w^t q(y,\theta)d\theta} q(z,w,dy\setminus \{z\}) \overline{P}^{(n-1)}(u,x;w,dz)dw,\  n = 1,2,\ldots\ .
\end{equation}

\begin{theorem}[\cite{FMS2}, Theorem 1]
\label{thm:JMP}
Given  a probability measure $\gamma$ on $\X$ and a $Q$-function $q$ satisfying Assumption~\ref{L1}, the jump process $\{\xi_t: t \in [T_0, T_1)\} $ defined in \eqref{jump} is a jump Markov process with the transition function $\overline{P}$.
\end{theorem}

\section{Kolmogorov's Equation's}\label{s3}
\textbf{3.1.} We start with \emph{Kolmogorov's backward equation:} for  $t\in(T_0,T_1),$  $x\in \X,$ and $B\in\B(\X),$
\begin{equation}
\label{eq:BKDE}
\frac{\partial}{\partial u}{P}(u,x;t,B) =
q(x,u){P}(u,x;t,B)
- \int_{\X }q(x,u,dy \setminus
\{x\})P(u,y;t,B) 
\end{equation} at $u\in [T_0,t).$ Since $q(x,u)$ is a real number, the right-hand side of \eqref{eq:BKDE} is always defined, and it is a real number if the function $P$ is bounded, as this takes place when $P$ is a transition function.

 Let $\cal P$ be the family of all real-valued non-negative functions $P(u,x;t,B),$ defined for all $t \in (T_0, T_1),$  $u \in [T_0, t),$ $x \in \X,$ and $B \in \B(\X),$  which are measurable  in $(u,x)\in [T_0, t)\times \X$ for all $t\in (T_0, T_1)$ and $B\in \B(\X).$ Observe  that $\overline{P} \in {\cal P}.$  

Consider a set $E$ and some family $\cal A$ of functions $f:E\to\overline{\BB{R}}=[-\infty,+\infty].$  A function $f$ from $\cal A$ is called minimal in the family $\cal A$ if for every function $g$ from $\cal A$ the inequality $f(x)\le g(x)$ holds  for all $x\in E.$

\textbf{3.2.} The following theorem describes the structure of solutions of Kolmogorov's backward equations.

\begin{theorem}[\cite{FMS2}, Theorem 2]
\label{thm:BKE}
Under Assumption~\ref{L1}, the transition function $\overline{P}$ is minimal in the family of the functions $P$ from $\cal P$ satisfying  the following two properties:

(i) for all $t\in(T_0,T_1),$  $x\in \X,$ and $B\in\B(\X),$
\begin{equation}
\label{BC1}
\lim\limits_{u \to t-}P(u,x;t,B) =  I \{ x \in B \},
\end{equation}
 and this function is absolutely continuous in $u \in [T_0, t);$

(ii) for each $t\in(T_0,T_1),$  $x\in \X,$ and $B\in\B(\X),$ Kolmogorov's backward equation \eqref{eq:BKDE}
holds at almost every $u\in [T_0,t).$

 In addition, if the transition function  $\overline{P}$ is  regular (that is, $\overline{P}(u,x;t,\X) $ $=1$ for all
$u,$ $x,$ $t$ in the domain of $\bar P$), then $\overline{P}$ is the unique
function in $\cal P$  satisfying  properties (i), (ii) 
and which is
 a measure on $(\X, \B(\X))$ for all $t\in(T_0,T_1),$ $u\in [T_0,t),$ and  $x\in \X,$
 and taking values in $[0,1]$.
 \end{theorem}

 We continue with \emph{Kolmogorov's forward equation:} for $u \in [T_0, T_1),$   $s \in (u,T_1),$  and $x \in \X,$ $B,$
\begin{equation}
\label{eq:FKDE}
\frac{\partial}{\partial t} P(u,x;t,B) = -\int_{B} q(y,t) P(u,x;t,dy)
 + \int_\X q(y,t,B\setminus \{y\})P(u,x;t,dy),
\end{equation}
at $t \in (u,s).$ If the jump rate $q(y,t)$ is bounded on $\X\times [T_0,T_1),$ then the right-hand side of \eqref{eq:FKDE} is defined, and it is finite if the function $P$ is bounded, as this take place for $P={\bar P}.$  If the function $q$ is not bounded, then the right-hand  side of  \eqref{eq:FKDE} can be undefined for a naturally defined function $P$  since  the right-hand  side of  \eqref{eq:FKDE}  is  $(-\infty)+(+\infty)$. The following example demonstrates this.
\begin{example}\label{Ex1}
{\rm \cite[Example 2]{FMS3}.(}$\X=\{0,\pm 1,\pm 2,\ldots\}$ and for $P=\overline{P}$ the right-hand side of Kolmogorov's forward equation~\eqref{eq:FKDE} is undefined for $B=\X,$ and for all $x\in\X$ it is finite when $B=\{z\}$ for each $z\in\X.${\rm ) Let us consider a homogeneous jump Markov chain, that is, the functions $q$ do not depend on the time parameter $t.$ We shall write $q(x):=q(x,t),$ $q(x,B):= q(x,t,B),$ and $q(x,j):=q(x,\{j\}).$  Let us denote $P(u,x;t,z):= P(u,x;t,\{z\}).$

We set $q(0) = 1,$  $q(0,  j)= 2^{-(|j| + 1)}$  and $q(j, -j) = q(j) =  2^{|j|}$ for all $ j \ne 0.$
 If $\xi_u=0,$ then starting at time $u$ the process spends at state 0 an exponentially distributed amount of time with the intensity $q(0)=1,$ then it jumps to a state $j\ne 0$ with probability $2^{-(|j|+1)},$ and then it oscillates between the states $j$ and $-j$ with equal intensities $2^{|j|}.$   Thus for all $u \in [T_0, T_1)$ and $t \in (u,T_1)$
\[\overline{P}(u,0;t,0) = e^{-(t-u)} \qquad \text{ and } \qquad \overline{P}(u,0;t,j)  = \frac{1- e^{-(t-u)}}{2^{|j|+1}}, \qquad j \ne 0,\]
where which implies 
\begin{multline*}
\int_\X q(y, \X \setminus \{y\})\overline{P}(u,0;t,dy)= \int_\X q(y)\overline{P}(u,0;t,dy) \\
= q(0)\overline{P}(u,0;t,0) + \sum_{j \ne 0} q(j)\overline{P}(u,0;t,j)= e^{-(t-u)} + \sum_{j > 0} (1- e^{-(t-u)}) = + \infty.
\end{multline*}
Thus, if $x=0$ and $B =\X$, then Kolmogorov's forward equation \eqref{eq:FKDE} does not not make sense with $P=\bar P$ because both integrals in \eqref{eq:FKDE} are infinite.  However, the integrals in \eqref{eq:FKDE} are defined for $B=\{z\},$ where $z\in \X.$ Indeed,  in this example
\[
\frac{\partial}{\partial t} P(u,0;t,0)=-P(u,0;t,0),
\]
 for $j\ne 0$
\[
\frac{\partial}{\partial t} P(u,0;t,j)=-2^{|j|+1}P(u,0;t,j)+2^{|j|}P(u,0;t,-j),
\]
and the right-hand side is a real number. It is easy to see that $\frac{\partial}{\partial t} P(u,i;t,j)$ is a real number for all $i,j\in\X.$\qed}
\end{example}

Example~\ref{Ex1} demonstrates that, when the right-hand side of Kolmogorov's forward equation \eqref{eq:FKDE} is not defined, this equation does not carry useful information, and this equation should be considered only for values of its parameters when its right-hand side  is well-defined. The following definition describes sets of states for which Kolmogorov's forward equation \eqref{eq:FKDE} is natural.

\begin{definition}
\label{defXnt}
For $s \in (T_0, T_1],$ a set $B\in \B(\X)$ is  called $(q,s)$-bounded if the function $q(x,t)$ is bounded on the set $B\times [T_0,s).$ 
\end{definition}

\begin{definition}
\label{defbbf}
A $(q,T_1)$-bounded set  is called $q$-bounded.
\end{definition}

Let $\hat{\cal P}$ be the family of real-valued functions $\hat{P}(u,x;t,B),$  defined for all $u \in [T_0, T_1)$, $t \in (u, T_1)$, $x \in \X,$ and $B \in \B(\X),$  which are measures on $(\X, \B(\ ))$ for fixed $u,$ $x,$ $t$ and are measurable functions in $t$ for fixed $u,$ $x,$ $B.$  In particular, $\bar P\in \hat{\cal P},$ where $\bar P$ is defined in \eqref{def}.

The following theorem plays the same role for Kolmogorov's forward equation as Theorem~\ref{thm:BKE} does for Kolmogorov's backward equation.
\begin{theorem}[{\rm \cite{FMS2}, Theorem 3}]
\label{thm:FKDE}
Under Assumption~\ref{ALB}, the transition function $\overline{P}$ is the minimal function in $\hat{\cal P}$ satisfying  the following two properties:

(i) for all $u \in [T_0, T_1),$   $s \in ]u,T_1),$  $x \in \X,$ and  $(q,s)$-bounded  sets $B,$
\begin{equation}
\label{BC2}
\lim_{t \to u+} P(u,x;t,B) = I\{x \in B\},
\end{equation}
and the function is absolutely continuous in $t \in ]u,s);$

(ii) for all $u \in [T_0, T_1),$   $s \in (u,T_1),$  $x \in \X,$ and  $(q,s)$-bounded  sets $B,$ Kolmogorov's forward equation \eqref{eq:FKDE}
holds for almost every $t \in (u,s).$

 In addition, if the transition function  $\overline{P}$ is  regular, then $\overline{P}$ is the unique
function in $\hat{\cal P}$  satisfying properties (i), (ii) 
and taking values in $[0,1]$. \end{theorem}

The following theorem provides the necessary and sufficient condition for a function $P\in\hat{\mathcal P}$ to satisfy Kolmogorov's forward equation.
\begin{theorem}[{\rm \cite{FMS2}, Theorem 4}]
\label{thm:intFKE}
 Let Assumption~\ref{ALB} hold. 
  A function $P$ from $\hat{\cal P}$ satisfies properties (i) and (ii)  stated in Theorem~\ref{thm:FKDE} if and only if, for all $u \in [T_0, T_1),$   $t \in (u,T_1),$  $x \in \X,$ and  $B \in \B(\X),$ 
\begin{equation}
\label{int-FKE}
\begin{split}
P(u,x;t,B) &= I\{x \in B\}e^{-\int_u^t q(x,\theta)d\theta} \\
&\qquad \qquad  + \int_u^t \int_\X \int_B e^{-\int_w^t q(y,\theta)d\theta} q(z,w, dy\setminus \{z\}) P(u,x;w,dz)  dw.
\end{split}
\end{equation}
\end{theorem}

Kolmogorov's forward equation can be also written in an integral form.
\begin{lemma}[{\rm \cite{FMS2}, Lemma 3}]
\label{Cor}
For arbitrary fixed $u \in [T_0, T_1)$, $s \in (u, T_1)$, $x \in \X$, and $(q,s)$-bounded set $B \in \B(\X),$ 
 a function $P$ from $\hat{\cal P}$ satisfies the  equality
\begin{equation}
\label{eq:FKE}
\begin{split}
P&(u,x;t,B)= I\{x \in B\} \\
&- \int_{u}^{t} \int_{B} q(y,w)P(u,x;w,dy)dw  + \int_u^t \int_\X q(y,w,B\setminus \{y\}) P(u,x; w, dy)  dw, \quad t \in (u,s),
\end{split}
\end{equation}
 if and only if it satisfies
the boundary  condition \eqref{BC2}, is absolutely continuous in $t \in (u,s),$ and satisfies Kolmogorov's forward equation~\eqref{eq:FKDE}  for almost every $t\in (u,s).$
\end{lemma}

Stronger results hold under Assumption~\ref{LB}, which is natural for CTJMDPs because it follows from Assumption~\ref{A1} introduced in the following section.  The following two statements describe necessary and sufficient conditions for the validity of Kolmogorov's forward equations under Assumption~\ref{LB}.

\begin{theorem} [{\rm \cite{FMS2}, Lemma 6}]
\label{LEMMA5}
Under  Assumption~\ref{LB}, a function  $P \in \hat{\cal P}$ satisfies properties (i) and (ii) stated in Theorem~\ref{thm:FKDE} if and only if
 if and only if the following two properties hold:

(a) for all $u \in [T_0, T_1)$, $x \in \X$, and  $q$-bounded sets $B$,   the  function $P(u,x;t,B)$ 
satisfies the boundary condition \eqref{BC2} and is absolutely continuous in $t \in (u,s)$ for each $s \in (u,T_1);$ 

(b) for all $u \in [T_0, T_1)$, $x \in \X$, and  $q$-bounded set $B$,   the function $P(u,x;t,B)$ satisfies Kolmogorov's forward equation~\eqref{eq:FKDE} for almost every $t\in (u,T_1).$
\end{theorem}
\begin{corollary}[{\rm \cite{FMS2}, Corollary 6}]
\label{Cor-S}
Under Assumption~\ref{LB}, the following statements hold:

(a) for all  $u \in [T_0, T_1),$ $x \in \X,$ and $q$-bounded sets $B \in \B(\X),$ the function $\overline{P}(u,x;t,B)$ satisfies  the equality in formula~\eqref{eq:FKE} for all $t \in (u,T_1).$ 

(b) the function $\overline{P}$ is the minimal function in $\hat{\cal P}$ for which statement~(a) holds. In addition, if the transition function $\overline{P}$ is regular, then $\overline{P}$ is the unique function in $\hat{\cal P}$ with values in $[0,1]$ for which statement~(a) holds.
\end{corollary}

Let us fix $x\in \X$  and $u = T_0.$ Then formula~\eqref{eq:FKE} becomes an equation with two variables $t$ and $B$. To simplify notations, we set  $P(t,B):=P(T_0,x;t,B)$  for any function $P$ from $\hat{\cal P}.$  Then \eqref{eq:FKE} becomes
\begin{equation}
\label{eq:FKE1}
P(t,B) = I\{x \in B\} + \int_{T_0}^{t} ds \int_\X q(y,s,B \setminus \{y\}) P(s, dy) -\int_{T_0}^{t} ds \int_{B} q(y,s) P(s,dy).
\end{equation}
For fixed $x \in \X$ and $u = T_0$, the function $\overline{P}(t, \cdot)$ is the marginal probability distribution 
of the  process $\{\xi_t: t \in [T_0, T_1)\}$ at time $t$ given $\xi_{T_0} = x$ and the initial state distribution is $\gamma.$ Under Assumption~\ref{LB}, the following corollary describes the minimal solution of \eqref{eq:FKE1} and provides a sufficient condition for its uniqueness.
\begin{corollary}[{\rm \cite{FMS2}, Corollary 6}]
\label{Cor-P}
Fix an arbitrary  $x \in \X$. Under Assumption~\ref{LB}, the following statements hold:

(a) for all $t \in (T_0, T_1)$ and $q$-bounded sets $B \in \B(\X),$ the function $\overline{P}(t,B)$ satisfies \eqref{eq:FKE1};

(b)  $\overline{P}(t,B),$ where $t\in (T_0,T_1)$ and $B\in \B(\X),$ is the minimal non-negative function $P(t,B),$ where $t \in (T_0, T_1)$ and $B\in\B(X),$  that is a measure on $(\X,\B(\X))$ for fixed $t$, is measurable in $t$ for fixed $B$, and satisfying \eqref{eq:FKE1} for all  $t \in (T_0, T_1)$ and for all $q$-bounded sets $B \in \B(\X).$ In addition, 
if $\overline{P}(t,\X) = 1$ for all $t \in (T_0, T_1)$,
then $\overline{P}(t,B)$  is the unique non-negative function  with values in $[0,1]$ and satisfying the conditions stated in the first sentence of this statement.
\end{corollary}
\section{Applications to Continuous-Time Jump Markov Decision Processes}. \label{s4}

\textbf{4.1.}  The probability structure of a CTJMDP  is specified by the four objects $\{ X, A, A(\cdot), \tilde{q} \}$, where
\begin{itemize}
\item[(i)] $(X, \B(X))$ is a standard Borel space (the state space);
\item[(ii)] $(A, \B(A))$ is a standard Borel space (the action space);
\item[(iii)] $A(x)$ is a non-empty subset of $A$ for each state $x \in X$ (the set of actions available at $x$). It is assumed that the set of feasible state-action pairs \[{\rm Gr}(A): = \{ (x,a) : x \in X, a \in A(x)\}\] is a measurable subset of $(X \times A)$ containing the graph of a measurable mapping of $ X$ to $A;$ 
\item[(iv)] $\tilde{q}( x, a, \cdot)$ is  a signed measure on $(X,\B(X))$ for each $(x,a)  \in {\rm Gr}(A)$ (the transition rate), such that $\tilde{q}(x,a, X) = 0$, $0 \le \tilde{q}(x,a, Z \setminus \{x\})  < \infty$, and $\tilde{q}(x,a,Z)$ is a measurable function on ${\rm Gr}(A)$ for each $Z \in \B(X).$
\end{itemize}

Let $\tilde{q}(x,a) := \tilde{q}( x,a,X \setminus \{x\} )$ for all $(x,a) \in {\rm Gr}(A),$ and let ${\bar q}(x):= \sup_{a \in A(x)}\tilde{q}(x,a)$ for all $x \in X$. If an action $a\in A(x)$ is selected at a state $x\in X$ and is fixed until the next jump, then the sojourn time has an exponential distribution with the intensity $\tilde{q}(x,a),$ and the process jumps to the set $Z\setminus\{x\},$  where $Z\in\B(X),$  with probability $\tilde{q}(x,a,Z\setminus\{x\})/\tilde{q}(x,a),$ if $\tilde{q}(x,a)>0,$ and  the state $x$ is absorbing if $\tilde{q}(x,a)=0.$  
In this paper we make everywhere the following standard assumption, which implies that there are no instantaneous jumps.
\begin{assumption}
\label{A1}
$\bar{q}(x) < \infty$ for each $x \in X$.
\end{assumption}

Similar to the case of jump Markov processes discussed in Section~\ref{s2}, we set $\bar{X}:= X\cup \{x_\infty\},$ where $x_\infty$ is an isolated point, and consider the sample space $(\Omega,\mathcal {F})$ with the filtration $\mathcal {F}_t,$ $t\ge 0,$ and stopping times $t_n$ representing jump epochs and $t_\infty=\lim_{n\to\infty}t_n,$ where the notation $\X$ is used in Section~\ref{s2} for the state space instead of $X.$  We also add an isolated point $a_\infty$ to the set of actions $A$ and set $\bar{A}:=A\cup\{a_\infty\}.$  $\mathcal{B}(\bar{A}):=\sigma(\mathcal{B}\{\bar{A}),\{a_\infty\}).$ Let $A(x_\infty):=\{a_\infty\}$ be the set of actions available at the state $x_\infty,$ and this set consists of the singleton $\{a_\infty\}.$  The state $x_\infty$ is always absorbing, that is, $\tilde{q}(x_\infty,a_\infty):=0.$  We recall that a mapping $p(\cdot|\cdot):E\times\mathcal{G}\to [0,1]$ is called a \emph{transition probability} from a measurable space $(E,\mathcal{E})$ to a measurable space $(G,\mathcal{G}),$ if $p(\cdot|e)$ is a probability on $(G,\mathcal{G})$ for all $e\in E,$ and $p(B|\cdot)$ is a measurable mapping of $(E,\mathcal{E})$ to $([0,1],\mathcal{B}([0,1])$ for every $B\in \mathcal{G}.$ Let $\P(\bar{A})$ be the set of probability measures in $(\bar{A},\mathcal{B}(\bar{A}).$

\textbf{4.2.} A policy $\pi$ is a mapping $\Omega\times \BB{R}_+\to  \P(\bar{A})$ such that:

(i) the stochastic process
  $\pi(B|\omega,t)$ is predictable for all $B\in\B(A),$ 

  (ii) $\pi(A( \xi_{t-}(\omega))|\omega,t)=1$ for all $(\omega,t)\in  \Omega\times\BB{R}_+$ with  $t< t_\infty(\omega),$ where the stochastic function $ \xi_{t}$ is defined in \eqref{jump} with $\X=X.$

   As follows from Jacod \cite[p. 241]{Jac}, the \emph{predictability} assumption means that there is a sequence of
transition probabilities $\pi^n:((X\times \BB{R}_+)^{n+1}, \B((X\times \BB{R}_+)^{n+1}))\to (A, \B(A))$ such that, at each  $t\in\BB{R}_+,$ the policy $\pi$ selects an action $a_t$ according to the probability measure
 \begin{equation}
\label{policy}
\pi(da_t| \omega, t) := \sum_{n \ge 0}\pi^n(da_t| x_0, t_1, x_1, \ldots, t_n, x_n, t-t_n)I\{ t_n < t \le t_{n+1}\} + \delta_{a_\infty} (da_t) I\{t \ge t_\infty\}, \ \  \omega \in \Omega,  
\end{equation}
where $\pi^n(A(x_n)| x_0, t_1, x_1, \ldots, t_n, x_n, t-t_n)=1 $ for $t\in (t_n,t_{n+1}],$ $n=0,1,\ldots,$  $\delta_{a_\infty}(\cdot)$ is a Dirac measure on $(\bar{A}, \B(\bar{A}))$ concentrated at $a_\infty,$ and we omit $\omega$ in the right-hand side of \eqref{policy} and in the condition following \eqref{policy}.  For example, the full version of the condition following \eqref{policy} is  $\pi^n(A(x_n(\omega))| x_0(\omega), t_1(\omega), x_1(\omega), \ldots, t_n(\omega), x_n(\omega), t-t_n(\omega))=1 $ for $t\in (t_n(\omega),t_{n+1}(\omega)],$ $n=0,1,\ldots.$ This condition  means that  $\pi(A( \xi_{t-}(\omega))|\omega,t)=1$ for all $(\omega,t)\in  (\Omega\times\BB{R}_+)$ with  $t< t_\infty(\omega).$

\textbf{4.3.}  A policy $\pi$ is called {\it Markov} if there exists  a  transition probability $\tilde{\pi}$ from $((X\times \BB{R}_+), \mathcal{B}((X\times \BB{R}_+)))$ to $(A, \B(A))$ such that
 $\pi(\cdot|\omega,t)=\tilde{\pi}(\cdot|\xi_{t-}(\omega),t)$ for all $(\omega,t)\in  (\Omega\times\BB{R}_+)$ with  $t< t_\infty(\omega)$. For a Markov policy $\pi$, formula \eqref{policy} implies that $\pi^n(B|x_0, t_1, x_1, \ldots, t_n, x_n, t-t_n)=\tilde{\pi}(B|x_n,t),$   when $t_n < t \le t_{n+1}$ and   for all $B\in\B(X)$ and  $n= 0,1,2,\ldots\ .$  With a slight abuse of notations, we shall write $\pi$ instead of $\tilde{\pi}.$

For $x\in\X,$ $p\in \P(A) $ such that $p(A(x)|x)=1,$ and $B\in \mathcal{B}(X),$ let us introduce the notations
\begin{equation*}
\label{q-def21}
\tilde{q}( x, p, B): =
\int_{A(z)} \tilde{q}(x, a, B) p (da).
\end{equation*}
and
\begin{equation}
\label{q-def22}
\tilde{q}(x,p): = \int_{A(z)}\tilde{q}(x,a) p(da);
\end{equation}
in general, for a measurable function $f$ on $X \times A$, we shall use the notation
\begin{equation}
\label{Ext}
f(x,p) := \int_{A(x)} f(x,a)p(da), \qquad x \in X,\ p \in \P(A),
\end{equation}
if the integral is defined.

Define the random measure $\nu^\pi$ on $(\BB{R}_+^0 \times X,\B(\BB{R}_+^0 \times X)),$ where $\BB{R}_+^0=(0,+\infty) ,$ by
\begin{equation}
\label{m2}
\nu^\pi (\omega; [ 0,t], B) := \int_{0}^{t} \tilde{q}( \xi_{s}(\omega), \pi_s(\omega), B \setminus \{\xi_{s}(\omega)\})I\{\xi_{s}(\omega) \in X\} ds,
\end{equation}
$\omega \in \Omega,$  $t \in \mathbf{R}_+,$ $B \in \B(X),$ where $\xi_s(\omega)$ is defined in \eqref{jump} with $T_0:=0$ and $T_1=+\infty,$ and $\pi_s(\omega)$ is the probability $\pi(\cdot|\omega,s),$ and we usually omit $\omega.$  This formula is similar to \eqref{compensator}. This random measure is predictable.  Indeed, in view of \eqref{jump} and \eqref{policy}, for each $Z\in\B(X),$ the stochastic process $\{\nu^\pi(\omega; [0,t], Z)\}$ is $\F_{t}$-measurable. In addition, it has continuous paths.  Therefore, these processes are $\F_{t-}$-measurable or, in other words, predictable; see, e.g., Jacod and Shiryaev~\cite[Proposition~2.6]{JS} or Kitaev and Rykov~\cite[Theorem 4.16]{KR}.

Furthermore, $\nu^\pi(\omega; [t_\infty, +\infty), X) = 0$ since $\xi_t(\omega) = x_\infty$ for all $t \ge t_\infty$ and $\nu^\pi(\omega; \{t\} \times X) = 0$ since the function $\nu^\pi(\omega; [ 0,t], X)$ is continuous in $t\in \mathbf{R}_+.$ In view of Jacod~\cite[Theorem~3.6]{Jac}, the predictable random measure $\nu^\pi$   and a probability measure $\gamma$ on $X$ define a unique probability measure ${\mathbf P}_\gamma^\pi$ on $(\Omega,\F)$ for which ${\mathbf P}_\gamma^\pi(dx_0)=\gamma(dx_0)$ and  $\nu^\pi$ is a compensator of the random measure  of the process $\xi_t.$  We remark that \cite[Theorem~3.6]{Jac} has two assumptions, namely, \cite[assumptions (4) and (A.2)]{Jac}.  Assumption (4) from \cite{Jac} is verified in the first sentence of this paragraph. Assumption (A.2) follows from the construction of the sample space $(\Omega,\F).$

\textbf{4.4.} Let us fix the initial state distribution $\gamma.$  This means that ${\mathbf P}_\gamma^\pi(\xi_0\in B)=\gamma(B)$ for all $B\in\mathcal{B}(X)$ and for each policy $\pi.$  Expectations with respect to the probability ${\mathbf P}_\gamma^\pi$ are denoted by $\mathbf{E}_\gamma^\pi.$  If the initial state distribution $\gamma$ is concentrated at a state $x\in X$, we shall write  ${\mathbf P}_x^\pi$ and $\mathbf{E}_x^\pi$ instead of   ${\mathbf P}_\gamma^\pi$ and ${\mathbf E}_\gamma^\pi$ respectively.

Kolmogorov's forward equation can be used to show that for many objective criteria for a given initial state distribution it is \emph{sufficient} to use a Markov policy.  This means that, if an initial state distribution is fixed, then for every policy $\pi$ there is a Markov policy $\varphi$ with the same or better value of the objective criterion.

We start with defining the Markov policy $\varphi.$  Let us define marginal distributions
\begin{equation}
\label{eqPszB}
P_\gamma^\pi(t,B):= {\mathbf P}_\gamma^\pi(\xi_t\in B)\qquad {\rm and}\qquad P_\gamma^\pi(t,B,U):= {\mathbf P}_\gamma^\pi(\xi_t\in B, a_t\in U),
\end{equation}
where $t\ge 0,$ $B\in \mathcal{B}(X),$ and $U\in\mathcal{B}(A).$ Of course, $P_\gamma^\pi(t,B)=P_\gamma^\pi(t,B,A).$  Since $P_\gamma^\pi(t,B,U)\le P_\gamma^\pi(t,B)$ for all $B\in \mathcal{B}(X),$  there is the Radon-Nikodym derivative $\frac{ P_\gamma^\pi(t,dx,U)}{P_\gamma^\pi(t,dx)}.$  Of course, as always, Radon-Nikodym derivatives $\frac{dP}{dQ}$ are defined $Q$-a.s.  The following lemma states that the defined Radon-Nikodym derivative can be presented by a Markov policy.  As follows from the definition provided above, a Markov policy $\phi$ is a transition probability $\phi(da|x,t)$ from $(\BB{R}_+\times X,\mathcal{B}(\BB{R}_+\times X))$ to $(A,\mathcal{B}(A))$ such that $\phi(A(x)|x,t)=1$ for all $x\in X$ and for all $t\in\BB{R}_+.$ The following lemma combines Lemma 1 and Corollary 1 from \cite{FMS3}.
\begin{lemma}\label{cordef}
For an initial state distribution $\gamma$  on $X$ and for a policy $\pi,$ there exists a Markov policy $\varphi$ such that for all   $t \in \BB{R}_+$  and $U \in \B(A),$ 
\begin{equation}
\label{formula-d}
\varphi(U \lvert x, t) = \frac{P_\gamma^{\pi}(t, dx,U)}{P_\gamma^{\pi}(t, dx)}, \quad  x \in X\  (P_\gamma^\pi(t, \cdot){\rm -}a.s.).
\end{equation}
\end{lemma}
\begin{proof} Since $\BB{R}_+\times X$ and $A$ are standard Borel spaces, and  $P_\gamma^\pi(t,B)=P_\gamma^\pi(t,B,A)$ for all $t\ge 0,$ and $B\in\mathcal{B}(X),$ in view of Bertsekas and Shreve~\cite[Corollary 7.27.1]{BS}, there exists a transition probability $\tilde{\varphi}$ from $(\BB{R}_+\times X,\B(\BB{R}_+\times X))$ to $(A,\B(A))$ such that formula~\eqref{formula-d} holds with $\varphi = \tilde{\varphi}.$ We observe that
\begin{equation}\label{eqAaas}
\int_X P_\gamma^\pi(t,dx)=P_\gamma^\pi(t,X)=\int_X\tilde{\varphi}(A|x,t) P_\gamma^\pi(t,dx)=\int_X\tilde{\varphi}(A(x)|x,t)P_\gamma^\pi(t,dx),
\end{equation}
where the second equality follows from the validity of \eqref{formula-d} for $\varphi=\tilde{\varphi},$
and the last equality follows from $\pi(A( \xi_{t-}(\omega))|\omega,t)=1$ for all $(\omega,t)\in  (\Omega\times\BB{R}_+)$ with  $t< t_\infty(\omega)$ and from \eqref{policy}.

For    a measurable mapping $\phi:X\to A$ with $\phi(x) \in A(x)$ for all $x \in \bar{X},$ whose existence is guaranteed by assumption (iii), let us define a Markov policy $\varphi,$
\begin{equation}
\label{varphi}
\varphi(U|x,t) := \left\{ \begin{array}{ll}
\tilde{\varphi}(U|x,t), \quad &\text{ if } \tilde{\varphi}(A(x)|x,t) =1,\\
\delta_{\phi(x)}(U), \quad &\text{ otherwise, }
\end{array}
\right.
\end{equation}
where $U\in \B(A)),$  $x \in X,$ and $t \in \BB{R}_+$.  Then \eqref{eqAaas} implies that $\tilde{\varphi}(A(x)|x,t)=1$, $x\in X$ ($P_\gamma^\pi(t, \cdot)$-a.s.), for all $t \in \BB{R}_+$. Therefore, $\varphi(U|x,t)=\tilde{\varphi}(U|x,t)$, $x\in X$ ($P_\gamma^\pi(t, \cdot)$-a.s.), for all $U\in\B(A)$ and $t \in \BB{R}_+$. Formula \eqref{formula-d} is proved.
\end{proof}

\textbf{4.5.} Kolmogorov's forward equation can be used to prove the following result.
\begin{theorem}[{\rm \cite{FMS3}, Theorem 1}]
\label{main}
For an initial distribution $\gamma$ on $X$ and a policy $\pi$, let $\varphi$ be a Markov policy satisfying \eqref{formula-d}. Then
\begin{equation}
\label{2-d}
P_\gamma^{\varphi} (t, B, U) \le P_\gamma^{\pi} (t, B, U), \qquad t \in \mathbf{R}_+, B \in \B(X), U \in \B(A).
\end{equation}
In addition, if $P_\gamma^\varphi(s,X) = 1$ for some $s \in \BB{R}_+$, then \eqref{2-d} holds for all $t\in (0,s]$ with an equality.  In particular, if $P_\gamma^\varphi(t,X) = 1$ for all $t \in \BB{R}_+$, then \eqref{2-d} holds  with an equality.
\end{theorem}
%
\begin{corollary}[{\rm \cite{FMS3}, Corollary 2}]
\label{C:main}
Let the transition rates $q(z,a)$ be bounded in  $(z,a)\in {\rm Gr}(A)$. Then, for every policy $\pi$ and initial distribution $\gamma$, $P_\gamma^\pi(t,X) = 1$ for all $t \in \BB{R}_+.$ In addition, formula~\eqref{2-d} holds with an equality for every Markov policy $\varphi$ satisfying \eqref{formula-d}.
\end{corollary}
\begin{corollary}[{\rm \cite{FMS3}, Corollary 3}]
For an initial distribution $\gamma$ on $X$ and a policy $\pi$, let $\varphi_1$ and $\varphi_2$ be two Markov policies satisfying \eqref{formula-d}. Then ${\mathbf P}_\gamma^{\varphi_1}={\mathbf P}_\gamma^{\varphi_2}$ and
\begin{equation*}
\label{2-deuf}
P_\gamma^{\varphi_1} (t, B, U) = P_\gamma^{\varphi_2} (t, Z, B), \qquad t \in \mathbf{R}_+,  B \in \B(X), U \in \B(A).
\end{equation*}
\end{corollary}
The answer to the question, whether inequality \eqref{2-d} holds in the form of an equality when $P_\gamma^\varphi(t,X)<1$, is open.

\textbf{4.6.} Here we explain the main ideas of the proof of Theorem~\ref{main}. The proof of Theorem~\ref{main} consists of two steps:

(i) proving the theorem for initial distributions concentrated at single points, that is $\gamma=\delta_x$ with $x\in X;$

(ii) extending the result obtained at step (i) to an arbitrary  distribution $\gamma$  of the initial state.

If the initial state $x\in X$ is fixed, every Markov policy $\phi$ defines the $Q$-function $q(y,s,B):= \tilde{q}(y,\varphi_s,B)$. where $y\in X,$ $s\in\BB{R}_+,$ and $B\in\B(X).$

Step (i) consists in proving the following statement.
\begin{lemma}[{\rm \cite{FMS3}, Lemma 2}]
\label{solutions}
For an initial state $x \in X$ and for a policy $\pi$, let $\varphi$ be a Markov policy satisfying \eqref{formula-d} with $\gamma(\{x\}) =1$.  Then, the following statements hold:


(i) for all $t\in\BB{R}_+$ and for all sets $B\in\B(X)$ such that $\sup_{\{z\in B,s\in \BB{R}_+\}} \tilde{q}(z,\phi_s)<+\infty,$  the functions $P(t,B)=P_x^\pi(t, B)$ and $P(t,B)=P_x^\varphi(t,B)$ satisfy Kolmogorov's forward equation
\begin{equation}
\label{MKDE}
P(t, B) = \delta_{x}(B) + \int_0^t  \int_X \tilde{q}(z,  \varphi_s, B\setminus \{z\}) P(s , dz)ds - \int_0^t \int_B \tilde{q}(z, \varphi_s)P(s,dz)ds;
\end{equation}
 (ii) for all $t \in \BB{R}_+$ and $B \in \B(X)$,
 \begin{equation}
 \label{1-d}
  P_x^\varphi(t, B) \le P_x^\pi(t, B);
 \end{equation}

 (iii) if $P_x^\varphi(s,X) = 1$ for some $s \in \BB{R}_+$, then \eqref{1-d} holds for $t \in (0,s]$ with an equality. In addition, if $P_x^\varphi(t,X) = 1$ for all $t \in \BB{R}_+$, then inequality~\eqref{1-d} holds with an equality for all $t \in \BB{R}_+$.
\end{lemma}

We would like to make the following two observations before providing the proof of Lemma~\ref{solutions}.  First, \eqref{MKDE} and \eqref{formula-d} imply \eqref{2-d}.  Second,
 $P(t,B)=P_x^\pi(t,B)$ satisfies \eqref{MKDE} for  $t\in\BB{R}_+$ and for  $B\in\B(X),$ such that $\sup_{\{z\in B,s\in\BB{R}_+\}}\tilde{q}(z,\varphi_s)<+\infty,$ because of the validity of formula \eqref{GKDE} stated in Lemma~\ref{l:GKE}. This formula is similar to Kolmogorov's forward equation, and  sometimes it is called Kolmogorov's forward equation.
For a policy $\pi$, a set $B \in \B(X)$ is called $(x, \pi)$-bounded if $\sup_{t\in \BB{R}_+} \mathbf{E}_x^\pi \tilde{q}(\xi_t, \pi_t) I\{\xi_t \in B\} < +\infty.$ 
\begin{lemma}[{\rm \cite{FMS3}, Lemma 4}]
\label{l:GKE}
For an initial state $x \in X$ and for  a policy $\pi$, the  formula
\begin{equation}
\label{GKDE}
P_x^\pi(t, B) = \delta_{x}(B)+ \mathbf{E}_x^\pi \int_{0}^{t} \tilde{q}(\xi_{s}, \pi_s, B\setminus \{\xi_s\})I\{\xi_s \in X\} ds - \mathbf{E}_x^\pi \int_{0}^{t} \tilde{q}(\xi_{s}, \pi_s)I\{\xi_s \in B\} ds
\end{equation}
 holds for all $t \in \BB{R}_+$ if the set $B\in \B(X)$ is $(x, \pi)$-bounded.
\end{lemma}

 Formula \eqref{GKDE} was introduced in \cite[Lemma 4]{Ki} for CTJMDPs with bounded rates $\tilde{q}(\cdot)$.  It was also proved in
 \cite[Theorem 3.1(c)]{GS} and in \cite[Theorem 1(b)]{PZ1} for non-explosive CTJMDPs satisfying different non-explosivity conditions and for $B\in\B(X)$ such that $\sup_{\{z\in B\}}\bar{q}(z)<+\infty;$ see also \cite[Proposition 4.29]{KR} and \cite[Theorem 2.4.5]{PZ}.

\begin{proof}[Proof of Lemma~\ref{solutions}]
We observe that equation \eqref{MKDE} is a particular case of equation \eqref{eq:FKE1} with $T_0=0,$ $T_1=+\infty,$ $\X=X,$ and $q(y,s,B):= \tilde{q}(y,\varphi_s,B)$ for $y\in X,$ $s\in\BB{R}_+,$ and $B\in\B(X).$ Since the Markov policy $\varphi$ defines jump Markov processes with a transition function $P(\cdot,\cdot,\cdot,\cdot)$ satisfying $P(x,0,t,B)=P_x^\varphi(t,B)$ for all $t\in\BB{R_+}$ and for all $B\in\B(X),$    we have that $\overline{P}(t,B')= P_x^\varphi(t,B')$ and, in view of Corollary~\ref{Cor-P}, $P_x^\varphi(t,B')$ is the minimal function $P(\cdot,\cdot)$ with the following properties: (a) the function $P(\cdot,B')$ is measurable for each $B'\in \B(X),$ (ii) $P(t,\cdot)$ is a measure on $(X,\B(X))$ for all $t\in\BB(R)_+,$ and (iii) equality \eqref{MKDE} holds for all  $t\in\BB(R)_+$ and for all $B\in \B(X)$ such that $\sup_{\{z\in B,s\in\BB{R}_+\}}\tilde{q}(z,\varphi_s)<+\infty.$

In order to complete the proof of Lemma~\ref{solutions}, we need to verify that \eqref{MKDE} holds  with $P(t,B)=P_x^\pi(t,B)$ for all $t\in\BB{R}_+$ and for all $B\in\B(X)$ such that $\sup_{\{z\in B,s\in\BB{R}_+\}}\tilde{q}(z,\varphi_s)<+\infty.$ If this is true, then statement (i) of the lemma is proved, it implies statement (ii) because $P_x^\varphi(t,B)$ is the minimal function described in the previous paragraph, and statement (iii) follows from Corollary~\ref{Cor-P}(b) applied to the $Q$-function $q$ defined in the previous paragraph.

For a non-negative measurable function $f$, for all $B \in \B(X)$ and $s \in \BB{R}_+$,
\begin{multline}
\label{f-int0}
\begin{aligned}
\mathbf{E}_x^\pi f(\xi_s, \pi_s) I\{\xi_s \in B\}  &= \int_B \int_{A(z)} f(z,a) P_x^\pi(s,dz,da)\\
&= \int_B \int_{A(z)} f(z,a) \varphi(da|z,s) P_x^\pi(s,dz) = \int_B  f(z,\varphi_s)  P_x^\pi(s,dz),
\end{aligned}
\end{multline}
where the first equality follows from \eqref{eqPszB}, the second equality follows from \eqref{formula-d}, and the last one follows from \eqref{Ext}. Then, for any non-negative measurable function $f$, for all $t \in \BB{R}_+$ and $B \in \B(X)$,
\begin{equation}
\label{f-int}
\mathbf{E}_x^\pi \left(\int_0^t f(\xi_s, \pi_s) I\{\xi_s \in B\} ds\right) =  \int_0^t  \mathbf{E}_x^\pi f(\xi_s, \pi_s) I\{\xi_s \in B\} ds = \int_0^t \int_B  f(z,\varphi_s)  P_x^\pi(s,dz) ds,
\end{equation}
where the first equality  follows from interchanging integration and expectation, and the second one follows from \eqref{f-int0}. Therefore, Lemma~\ref{l:GKE}, formula \eqref{f-int} with $B = X$ and $f(\xi_s, \pi_s) = \tilde{q}(\xi_s, \pi_s, B \setminus \{\xi_s\}),$  and the same formula with $f(\xi_s, \pi_s) = \tilde{q}(\xi_s, \pi_s)$ imply that the function $P_x^\pi(t,B)$ satisfies Kolmogorov's forward equation~\eqref{MKDE}  for all $t \in \BB{R}_+$ if the set $B \in \B(X)$ is  $(x,\pi)$-bounded.

To conclude the proof of the lemma, we need to check that, if $\sup_{\{z \in B, s \in \BB{R}_+\}}\tilde{q}(z,\varphi_s) <  +\infty$ for  $B\in\B(X),$ then the set $B$ is $(x,\pi)$-bounded.   This is true because, if $B\in\B(X)$ and $\sup_{\{z \in B, s \in \BB{R}_+\}}\tilde{q}(z,\varphi_s) <  +\infty,$ then
\begin{multline*}
\label{bound}
\sup_{s \in \BB{R}_+}\mathbf{E}_x^\pi  \tilde{q}(\xi_s, \pi_s) I\{\xi_s \in B\} =  \sup_{s \in \BB{R}_+} \int_B \tilde{q}(z,\varphi_s) P_x^\pi(s,dz) \le \left(\sup_{z \in B, s \in \BB{R}_+} \tilde{q}(z,\varphi_s) \right)P_x^\pi (s, B) <  +\infty,
\end{multline*}
where the first equality follows from \eqref{f-int0} with $f(\xi_s, \pi_s) = \tilde{q}(\xi_s, \pi_s)$, the first inequality is straightforward, and the last one is true since $\sup_{\{z \in B, s \in \BB{R}_+\}}\tilde{q}(z,\varphi_s) <  +\infty,$ and $P_x^\pi(s, B) \le 1$. Therefore, the function $P_x^\pi(t,B)$ is a solution of Kolmogorov's forward equation~\eqref{MKDE}, and statement (i) of the lemma holds.
\end{proof}

\textbf{4.7.} For step (ii) of the proof of Theorem~\ref{main}, which is proving Theorem~\ref{main} for an initial state distribution $\gamma,$ if it is known that the theorem is valid for any initial state $x,$ we expand the state and the action sets $X$ by adding an additional point $x'\notin X$ to the state space $X$ and by adding two points $a'\notin A$ and $a''\notin A$ to the action set $A.$  So, $X':=X\cup\{x'\}$ and $A':=A\cup\{a',a''\}$ are the expanded state and action sets.   In addition, the set of actions at the added state $x'$ is $A'(x'):=\{a',a''\},$ and the action  $a''$ is added to each action set $A(x)$ for $x\in X.$  In other words, the new action sets at $x\in\X$ are $A'(x):=A(x)\cup\{a''\}.$

If the action $a'$ is permanently chosen at state $x',$ then the process jumps with intensity 1 to $X$, and the distribution of the next state is $\gamma.$  If the action $a''$ is permanently chosen at a state $x\in X',$ then the state $x$ is absorbing.  All other actions behave is in the same way as for the original model.  So, for the given function $\tilde{q}$ for the original model, we define the function $q'$ for the expanded model: for all $x \in X', a \in A'(x),$ and $B \in \B(X'),$
\begin{equation}
\label{new-q1}
q'( x,a,B) := \left \{
\begin{array}{ll}
\tilde{q}( x, a, Z \setminus \{x'\} ), &\quad \text{ if } x \in X, a \in A(x),\\
\gamma(B\setminus \{x'\})-\delta_{x'}(B)   &\quad \text{ if }  x = x', a = a',\\
0, &\quad \text{ if } x \in X', a = a''.
\end{array}
\right.
\end{equation}

Let us choose an arbitrary policy $\pi$ and an initial state distribution $\gamma$ for the original model.  Let us fix an arbitrary constant $u>0.$  We shall construct a special policy $\tilde{\pi}$ for the new model.  For our purposes, it is sufficient to define  the policy $\gamma$ only for the initial state $x'.$ At the state $x'$, the policy $\tilde{\pi}$  chooses the action $a'$ at time $t,$ if  $t<u,$ and $\tilde{\pi}$  chooses action the action $a''$ if time $t,$ if $t\ge u.$ At each state $x\in X,$ the policy  $\tilde{\pi}$ chooses the action $a''$ if $t<u,$ and it observes the process $\xi_t$ starting from time $u$ until $t,$ if $t>u$ and $x_t\in X,$ and chooses actions with the same probabilities as the policy $\pi$ would choose at time $(t-u)$ on the basis of the same observations.  The formal definition of the policy $\tilde{\pi}$ is provided in \cite{FMS}.

Starting from the state $x'$, under the policy $\tilde{\pi}$ the process either jumps from the state $x'$ during the time interval $[0,u)$ or stays for good at the state $x'$ with probability 1.  Of course, the probability that the process jumps at $t=u$ is 0.  If the process jumps from $x'$ to a state $x\in\X$ at time $t<u,$ then the policy $\tilde{\pi}$ keeps the process at the state $x$ until time $t=u,$ and then resets the clock to time 0 and starting from time $t=u$ the policy $\tilde{\pi}$ behaves in the same way as the policy $\pi$ staring from time 0. Therefore, the process always stays at state $x'$ with the probability $e^{-u},$ and $P_{x'}^{\tilde{\pi}}(\xi_{t+u}\in B,a_{t+u}\in U)=(1-e^{-u})P_\gamma^\pi(\xi_t\in B,a_t\in  U)$ for $B\in\B(X)$ and $U\in \B(A).$  Thus, the measures $P_{x'}^{\tilde{\pi}}((\xi_t,a_t)\in D)$ and  $P_\gamma^\pi((\xi_t,a_t)\in D),$ where $D\in\B(X\times A),$ are equivalent.

Let $\varphi$ be a Markov policy satisfying \eqref{formula-d}.  Then the Markov policy $\tilde{\varphi}$ for the new model satisfies
$\tilde{\varphi}(B|x,t+u)=\frac{P_\gamma^{\varphi}(t,dx, B)}{P_\gamma^{\varphi}(t,dx)}=\varphi(B|x,t)$ for $t\ge 0,$ and it may be chosen selecting the action $a'$ at state $x'$ and action $a''$ at states $x\in X$ for $t<u.$ Thus, we also have $P_{x'}^{\tilde{\varphi}}(\xi_{t+u}\in B,a_{t+u}\in U)=(1-e^{-u})P_\gamma^\varphi(\xi_t\in B,a_t\in  U)$ for $B\in\B(X)$ and $U\in \B(A).$ We have that
\begin{equation}
\label{2-ddddd}
P_\gamma^{\varphi} (t,B) \le P_\gamma^{\pi} (t, B), \qquad t \in \mathbf{R}_+, B \in \B(X)
\end{equation}
 since $P_{x'}^{\tilde{\varphi}}(x_{t+u}\in B)\le P_{x'}^{\tilde{\pi}}(x_{t+u}\in B).$   The latter equality holds because of the validity of Theorem~\ref{main} for an initial state distributions concentrated at a single point $x'.$ Formula \eqref{2-ddddd} implies \eqref{2-d} since for $t\in\BB{R}_+,$ $B\in\B(X),$ and $U\in \B(A),$
\[
P_\gamma^\varphi(t,B,U)=\int_B \varphi(U|y,t)P_\gamma^\varphi(t,dy)=\int_B \frac{P_\gamma^\pi(t,dy,U)}{P_\gamma^\pi(t,dy)}P_\gamma^\varphi(t,dy)\le \int_B P_\gamma^\pi(t,dy,U) =P_\gamma^\pi(t,B,U),
\]
where the first inequality follows from the definition of a Markov policy, the second equality follows from \eqref{formula-d}, the inequality follows from  \eqref{2-ddddd} and the Radon-Nikodym theorem, and the last equality is obvious.

To prove the last statement of Theorem~\ref{main}, assume that $P_\gamma^\varphi(s,X) = 1$ for some $s \in \BB{R}_+$. We fix an arbitrary $t\in (0,s].$ Then $P_\gamma^\varphi(t,X) = 1.$

Let $P_\gamma^\varphi(t,B,U)<P_\gamma^\pi(t,B,U)$ for some $B\in\B(X)$ and $U\in\B(A).$ Then, in view of \eqref{2-d}, $P_\gamma^\varphi(t,B,A\setminus U)\le P_\gamma^\pi(t,B,A\setminus U)$ and
$P_\gamma^\varphi(t,X\setminus B,A)\le P_\gamma^\pi(t,X\setminus B,A).$  Therefore
\[\begin{aligned}
&1= P_\gamma^\varphi(t,X)=P_\gamma^\varphi(t,B,U)+P_\gamma^\varphi(t,B,A\setminus U)+P_\gamma^\varphi(t,X\setminus B,A)\\ &<P_\gamma^\pi(t,B,U)+P_\gamma^\pi(t,B,A\setminus U)+P_\gamma^\pi(t,X\setminus B,A)=P_\gamma^\pi(t,X),
\end{aligned}
\]
which is impossible since $P_\gamma^\pi(t,X)\le 1$.  Thus,   $P_\gamma^\varphi(t,B,U)=P_\gamma^\pi(t,B,U)$ for all $B\in\B(X)$ and for all $U\in\B(A).$

\section{Applications of Theorem~\ref{main} to Particular Objective Criteria} \label{s5}
\textbf{5.1.} For a cost function $c:X\times A\to (\BB{R}_+,\B(\BB{R}_+)),$ let us consider the expected total costs over the infinite time horizon
\begin{equation}\label{einitedc}
V_\alpha(\gamma,\pi):=\E_\gamma^\pi\int_0^{t_\infty} e^{-\alpha t}c(\xi_t,\pi_t)dt,
\end{equation}
where the constant $\alpha>0$ is the discount rate, $\gamma$ is the initial state distribution, $\pi$ is a policy, and $c(\xi_t,\pi_t)$ is the cost rate at time $t,$ where $c(\xi_t,\pi_t)$ is defined by formula \eqref{Ext} with $c=f,$ $x=\xi_t,$ and $p=\pi_t.$

Since the state  $x_\infty$ is always absorbing, we can rewrite $V_\alpha(\gamma,\pi)$ as
\begin{equation}\label{einitdp}
V_\alpha(\gamma,\pi):=\E_\gamma^\pi\int_0^{+\infty} e^{-\alpha t}c(\xi_t,\pi_t)I\{\xi_t\in X\}dt=\int_0^{+\infty}e^{-\alpha t} \E_\gamma^\pi[c(\xi_t,\pi_t)I\{\xi_t\in X\}]dt,
\end{equation}
where in the last equality we exchanged the order of the expectation and integration.  Since
\[\E_\gamma^\pi[c(\xi_t,\pi_t)I\{\xi_t\in X\}]=\int_X\int_a c(z,a)P_\gamma^\pi(t,dz,da),
\]
 for a Markov policy $\varphi$ satisfying \eqref{formula-d}, formula \eqref{einitdp} implies
\[
V_\alpha(\gamma,\pi)=\int_0^{+\infty}e^{-\alpha t}\left [ \int_X\int_a c(z,a)P_\gamma^\pi(t,dz,da)\right ]dt\] \[\ge \int_0^{+\infty}e^{-\alpha t}\left [ \int_X\int_a c(z,a)P_\gamma^\varphi(t,dz,da)\right ]dt=V_\alpha(\gamma,\varphi),
\]
where the inequality follows from Theorem~\ref{main}.  In addition, if $P_\gamma^\varphi(t,X)=1$ for all $t>0,$ then the inequality becomes an equality and $V_\alpha(\gamma,\pi)=V_\alpha(\gamma,\varphi).$

Thus,  Theorem~\ref{main} implies that the total expected cost defined in \eqref{einitedc} is smaller or equal for a Markov policy satisfying \eqref{formula-d} than for the original policy $\pi$ if the initial state distribution $\gamma$ is fixed.  The same conclusions hold for some other objective criteria. In the rest of this section we discuss two additional objective criteria, and one of them is more general than \eqref{einitedc}.

\textbf{5.2.} Let the time horizon is finite, that is $t\in [0,T]$ with $T<+\infty.$  In addition to the cost rates $c$ described above, the costs also are incurred at deterministic time instances $(u_i\in[0,T])_{i=1,2,\ldots}.$  The instant costs are defined by the measurable functions $G_i:(X\times A)\to (\BB{R}_+,\B(\BB{R}_+),$ $i=1,2,\ldots$.  A variable discount rate $\alpha(\cdot):[0,t]\to \BB{R}$ is a Borel function. We always assume that $\int_0^T\alpha(t) dt$ is defined, that is, either  $\int_0^T\alpha^+(t) dt<+\infty$ or $\int_0^T\alpha^-(t) dt>-\infty,$ where $d^+:=\max\{d,0\}$ and $d^-:=\min\{d,0\}$ for a number $d.$

In this case, for an initial state distribution $\gamma$ and a policy $\pi,$ the \textit{finite-horizon expected total discounted cost} up to time $T \in \BB{R}_+$  with a variable discount rate $\alpha(\cdot)$ is
\begin{equation}
\label{FDR}
V_{\alpha}^{T}(\gamma, \pi) :=  \mathbf{E}_\gamma^\pi \left[\int_{0}^{T \wedge t_\infty} e^{-\int_0^t\alpha(s) ds}c(\xi_{t}, \pi_t)dt + \sum\limits_{i=1}^{\infty} e^{-\int_0^{u_i}\alpha(s) ds}G_i(\xi_{u_i}, \pi_{u_i})I\{\xi_{u_i}\in X\}\right].
\end{equation}
In particular, if $u_1=T$ and $G_i\equiv 0$ for all $i>1,$ we deal with the problem with the terminal cost $G_1$ collected at the termination time $T.$ Because of the same arguments, we also have that $V_\alpha^T(\gamma,\varphi)\le V_\alpha^T(\gamma,\varphi)$ for a Markov policy $\varphi$ satisfying \eqref{formula-d}, and the equality takes place if $P_\gamma^\varphi(T,X)=1.$

\textbf{5.3.} For the infinite-horizon $T=+\infty,$ a more general objective criterion than the one defined in \eqref{einitedc} can be considered.  We consider the nonnegative cost functions $c$ and $G_i$ defined above, where costs $G_i,$ $i=1,2,\ldots,$ are collected at deterministic time instances $u_i\in \BB{R}_+.$  In addition, costs $C(\xi_{t_{n-1}}, \xi_{t_n})$ are collected at jump epochs $t_n,$ where $n=1,2,\ldots,$ where $C$ is a nonnegative Borel function $C:X\times X\to \BB{R}_+.$  The discount rate $\alpha>0$ is constant.  The \emph{infinite-horizon expected total discounted cost} is
\begin{equation*}
\label{IDR}
V_{\alpha}(\gamma, \pi): =  \mathbf{E}_\gamma^{\pi} \left[\int_{0}^{t_\infty} e^{-\alpha t}c(\xi_{t}, \pi_t)dt + \sum\limits_{n=1}^\infty e^{-\alpha t_n}C(\xi_{t_{n-1}}, \xi_{t_n}) + \sum\limits_{i=1}^{\infty} e^{-\alpha {u_i}}G_i(\xi_{u_i}, \pi_{u_i})I\{\xi_{u_i}\in X\}\right].
\end{equation*}
In this case it is also true that $V_\alpha(\gamma,\varphi)\le V_\alpha(\gamma,\varphi)$ for a Markov policy $\varphi$ satisfying \eqref{formula-d}, and the equality takes place if $P_\gamma^\varphi(t,X)=1$ for all $t\in\BB{R}_+;$ see \cite[Theorem 5]{FMS3}. As shown in \cite[Example 1]{FMS3}, this may not be true if the cost function $C$ also depends on an action chosen at jump epochs. Other criteria and additional results can be found in \cite[Section 6]{FMS3}.




\begin{thebibliography}{99}
\bibitem{Anderson}
W. J. Anderson,   {\it Continuous-Time Markov Chains: An Applications-Oriented Approach}.  Springer-Verlag, New York, 1991.
\bibitem{BS}
D. P. Bertsekas, S. E.  Shreve,   {\it Stochastic Optimal Control: The Discrete-time Case}.   Academic Press, New York, 1978.

\bibitem{DS}
C. Derman, R. Strauch, \emph{A note on memoryless rules for controlling
  sequential control processes,} Ann. Math. Statist., 37(1966), pp.276-278.

\bibitem{Do}
J.L. Doob, \emph{Markoff chains~--~denumerable case,} Trans. Amer. Math. Soc. 58 (1945), pp. 455-473.

\bibitem{Fe}
E. A. Feinberg,  \emph{A generalization of “expectation equals reciprocal of intensity” to non-stationary
exponential distributions,} J. Appl. Probab., 31 (1994), pp. 262-267.

\bibitem{Fe1}
E. A.  Feinberg,  \emph{Continuous time discounted jump Markov decision processes: a discrete-event
approach,} Math. Oper. Res. 29 (2004), pp. 492-524.

\bibitem{Fe2}
E.A. Feinberg,  \emph{Reduction of discounted continuous-time mdps with unbounded jump and reward
rates to discrete-time total-reward MDPs.} In D. Hernandez, A. Minjarez (Eds.), \emph{Optimization, Control, and
Applications of Stochastic Systems,}   pp. 201-213, 2012, Birkhäuser/Springer, New York.

\bibitem{FMS}
  E. A. Feinberg, M. Mandava, A. N. Shiryaev,  \emph{On solutions of Kolmogorov's equations for nonhomogeneous jump Markov processes,}  J. Math. Anal. Appl., 411 (2014), pp. 261-270.

 \bibitem{FMS2}
  E.A. Feinberg, M. Mandava, A.N. Shiryaev, \emph{ Kolmogorov’s Equations for jump Markov        processes with unbounded jump rates,} Ann. of Oper. Res., (2017), published       online: DOI 10.1007/s10479-017-2538-8.

 \bibitem{FMS3}
  E.A. Feinberg, M. Mandava, A.N. Shiryaev,  \emph{Sufficiency of Markov policies for continuous-time jump Markov decision processes,} \emph{Math. Oper. Res.,} (2021), published online:
  https://doi.org/10.1287/moor.2021.1169.


\bibitem{Fel}
W. Feller,   \emph{On the integro-differential equations of
purely-discontinuous Markoff processes,}  Tran. Amer. Math.
Soc., 48 (1940),  pp. 488-515; \emph{Errata,}  Trans. Amer.
Math. Soc., 58 (1945), p. 474.

\bibitem{GHL}
X. Guo,  O. Hern\'{a}ndez-Lerma, \emph{Continuous-time controlled {M}arkov chains,} Ann. Appl. Probab., 13(2003), pp. 363-388.

\bibitem{GH}
 X. Guo,  O. Hern\'{a}ndez-Lerma,   {\it Continuous-Time Markov Decision
Processes: Theory and Applications.}  Springer-Verlag, Berlin, 2009.

\bibitem{GR}
X. Guo, U. Rieder, \emph{Average optimality for continuous-time controlled {M}arkov processes in {P}olish spaces,} Ann. Appl. Probab., 16(2006), pp. 730-756.

\bibitem{GS}
X. Guo, X. Song, \emph{Discounted continuous-time constrained {M}arkov decision
  processes in {P}olish spaces,} Ann. Appl. Probab., 21(2011), pp. 2016-2049.

  \bibitem{GZ}
Xin Guo, Y. Zhang, \emph{A useful technique for piecewise deterministic Markov decision processes,} Oper. Res. Let., 49(2021), p. 66-61.

\bibitem{Ho}
R.A. Howard, \emph{Dynamic Programming and Markov Processes,} Wiley, New York, 1960,

\bibitem{Jac}
J. Jacod,  \emph{ Multivariate point processes: predictable projection, Radon-Nikodym derivatives, representation of martingales,}  Probab. Theory Related Fields, 31 (1975), pp. 235-253.

\bibitem{JS}
 J. Jacod, A. N. Shiryaev, \emph{Limit Theorems for Stochastic Processes}. Springer-Verlag, New York, 2003.

%

\bibitem{Ka}
P. Kakumanu, \emph{Continuously discounted Markov decision model with countable state and action
space,} Ann. Math. Stat., 42(1971), pp. 919-926.

\bibitem{Kechris}
A. S. Kechris,   {\it Classical Descriptive Set Theory}. Springer, New York, 1995.

\bibitem{Ken}
D. G. Kendall,    \emph{Some further pathological examples in
the theory of denumerable Markov processes,}  Q. J. Math., 7 (1956), pp. 39-56.



\bibitem{Ki}
M. Yu. Kitaev. \emph{Semi-Markov and jump Markov controlled models: average cost
  criterion,} Theory Prob. Appl., 30((1985)), pp. 272-288.


\bibitem{KR}
 M. Yu. Kitaev, V. V. Rykov, \emph{Controlled Queueing Systems.} CRC Press, Boca Raton, 1995.

\bibitem{Kol}
A. N. Kolmogorov,    \emph{On analytic methods in probability theory} (in German 1931, in
Russian 1938). In  A.N. Shiryaev (Ed.), {\it Selected Works of A.N.
Kolmogorov}, Vol. II, Probability Theory and Mathematical
Statistics. Springer, New York, pp. 62-108, 1992.


\bibitem{Kuz}
S. E. Kuznetsov,   \emph{Any Markov process in a Borel space
has a transition function.}   Theory Probab. Appl., 25 (1981), pp. 384-388.

\bibitem{Mi}
B. Miller, \emph{Finite state continuous time Markov decision processes with a finite planning horizon,}
SIAM J. Control, 6(1968), pp. 266-280.

\bibitem{Mi1}
B. Miller, \emph{Finite state continuous time Markov decision processes with an infinite planning horizon,}
J. Math. Anal. Appl., 22(1968), pp. 552-569.


\bibitem{Pi}
A.B. Piunovskiy, \emph{Realizable strategies in continuous time Markov decision processes,} SIAM J. Control Optim., 56(2018), pp. 473-495.

\bibitem{PZ1}
A.B. Piunovskiy, Y. Zhang, \emph{Discounted continuous-time {M}arkov decision
  processes with unbounded rates and randomized history-dependent policies: the
  dynamic programming approach.} 4OR-Q. J. Oper. Res., 12(2014), pp. 49-75.

  \bibitem{PZ}
A.B. Piunovskiy, Y. Zhang, \emph{Continuous-Time Markov Decision Processes.} Springer Nature, Switzerland, 2020.

\bibitem{Reu}
G. E. H.    Reuter,   \emph{Denumerable Markov processes and the
associated contraction semigroups on $l$.}  Acta Math., 97 (1957), pp. 1-46.

\bibitem{Ry}
V.V. Rykov, \emph{Markov decision processes with finite state and decision spaces,} Theory Probab. Appl., 11 (1966), pp. 302-311.

\bibitem{Se}
R. F.  Serfozo,  \emph{An equivalence between continuous and discrete time Markov decision processes,} Oper. Res.
27(1979), pp. 616-620.


\bibitem{St}
R. Strauch,\emph{ Negative dynamic programming,} Ann. Math. Statist.,
  37(1966), pp. 871--890.




\bibitem{YGHL}
 L. Ye, X. Guo,  O. Hern\'{a}ndez-Lerma,   \emph{Exstence and regularity of a nonhomogeneous transition matrix under measurability conditions.} J. Theoret. Probab., 21 (2008), pp. 604-627.

 \bibitem{Yus}
A. Yushkevich, \emph{Controlled {M}arkov models with countable state space and
  continuous time,} Theory Probab. Appl., 22(1977), pp. 215-235.

 \bibitem{Yus1}
A. Yushkevich,  \emph{Controlled jump {M}arkov models,} Theory
  Probab. Appl., 25(1980), pp. 244-266.



\end{thebibliography}
\end{document}